\documentclass[10pt]{article}

\usepackage{amsmath,amssymb,latexsym}
\usepackage{amsthm}
\usepackage{array,graphicx,mathrsfs}
\usepackage{hyperref,stmaryrd,xcolor} 
\usepackage{dsfont}

\newtheorem{theorem}{Theorem}[section]

\newtheorem{corollary}[theorem]{Corollary}  
\newtheorem{proposition}[theorem]{Proposition}
\newtheorem{lemma}[theorem]{Lemma}
\newtheorem{remark}[theorem]{Remark}

\numberwithin{equation}{section}

\newcommand{\ds}{\displaystyle}

\newcommand{\ld}{{L_2}}

\newcommand{\norm}[1]{\left\Vert#1\right\Vert}
\newcommand{\ic}{{\bf i}}

\title{On the reachable set for the one-dimensional heat equation}
\author{
	J\'er\'emi Dard\'e\footnote{e-mail: {\tt jeremi.darde@math.univ-toulouse.fr}.} 
	\qquad  
	Sylvain Ervedoza\footnote{e-mail: {\tt sylvain.ervedoza@math.univ-toulouse.fr}.}
	\smallskip
	\\
	{\it\footnotesize Institut de Math\'ematiques de Toulouse ; UMR 5219 ;  Universit\'e de Toulouse ; CNRS ;}
	\\
	{\it\footnotesize UPS IMT F-31062 Toulouse Cedex 9, France}
}
\date{\today}
\begin{document}
\maketitle
\abstract{
	The goal of this article is to provide a description of the reachable set of the one-dimensional heat equation, set on the spatial domain $x \in (-L,L)$ with Dirichlet boundary controls acting at both boundaries. Namely, in that case, we shall prove that for any $L_0 >L$ any function which can be extended analytically on the square $\{ x + \ic y,\, |x| + |y| \leq L_0\}$ belongs to the reachable set. This result is nearly sharp as one can prove that any function which belongs to the reachable set can be extended analytically on the square $\{ x + \ic y,\, |x| + |y| < L\}$. Our method is based on a Carleman type estimate and on Cauchy's formula for holomorphic functions.
}
%
%
\section{Introduction}
{\bf Setting.} The goal of this article is to describe the reachable set for the $1$d heat equation. To fix the ideas, let $L,\, T >0$ and consider the equation
\begin{equation}
	\label{Eq-Heat-1d}
		\left\{
			\begin{array}{ll}
				\partial_t u - \partial_{xx} u = 0 & \hbox{ in } (0,T) \times (-L,L), 
				\\
				u(t,-L) = v_-(t) & \hbox{ in } (0,T), 
				\\
				u(t,L) = v_+(t)  & \hbox{ in } (0,T), 
				\\
				u(0, x) = 0 & \hbox{ in } (-L,L).
			\end{array}
		\right.
\end{equation}
In \eqref{Eq-Heat-1d}, the state $u = u(t,x)$ satisfies a heat equation controlled from the boundary $x \in \{-L, L\}$ through the control functions $v_-(t),\, v_+(t)  \in L^2(0,T)$. In this article, the control functions $v_-$ and $v_+$ will be complex valued unless stated otherwise, and following, the solutions of the heat equation will also be complex valued.
\\
Our goal is to describe the reachable set $\mathscr{R}_L(T)$ at time $T>0$, defined as follows:
\begin{equation}
	\label{Def-Reachable-T}
		\mathscr{R}_L(T)
		=
		\{u(T)\ |\  u \hbox{ solving \eqref{Eq-Heat-1d} with control functions } v_-, \, v_+ \in L^2(0,T) \}.
\end{equation}
Obviously, due to the linearity of the problem, the reachable set $\mathscr{R}_L(T)$ is a vector space. Besides, as remarked by Seidman in \cite{Seidman}, due to the fact that the heat equation is null-controllable in arbitrarily small time (\cite{FattoriniRussel71} in dimension one, \cite{LebRob,FursikovImanuvilov} in higher dimensions), the reachable set $\mathscr{R}(T)$ does not depend on the time horizon $T$. Using again the null-controllability of the heat equation in small time, one also easily checks that the set of states $u(T, \cdot)$ which can be reached by solutions of \eqref{Eq-Heat-1d}$_{(1,2,3)}$ starting from an initial datum $u(0, \cdot) \in L^2(-L,L)$ coincides with $\mathscr{R}$.
Therefore, we will simply denote the reachable set by $\mathscr{R}$ in the following.
\\�\par
\noindent{\bf Main results.} We aim at proving the following result, whose proof is given in Section \ref{Sec-Proof-Main}:
\begin{theorem}
	\label{Thm-Main}
		For $L_0 > 0$, let us introduce the (open) square
		\begin{equation}
			\label{Def-Square-R}
			\mathcal{S}(L_0) = \{ x+ \ic y, \,  |x| + |y| < L_0\}, 
		\end{equation}
		and let us denote by $\mathscr{A}(L_0)$ the set of functions $u \in L^2(-L,L)$ which can be extended analytically to the set $\mathcal{S}(L_0)$:
		\begin{equation}
			\label{Def-A-R}
			\mathscr{A}(L_0) = \{ u \in L^2(-L,L) \hbox{ which can be extended analytically to } \mathcal{S}(L_0)\}.
		\end{equation}
		Then we have the following:
		\begin{equation}
			\label{Subset-Reachable}
			\underset{L_0 >L}{\cup} \mathscr{A}(L_0) \subset \mathscr{R}_L .
		\end{equation}
\end{theorem}
Theorem \ref{Thm-Main} describes a subset of the reachable set $\mathscr{R}$ in terms of  analytic extensions of functions. It turns out that using this description, Theorem \ref{Thm-Main} is \emph{sharp} as a consequence of the following result, recently obtained in \cite[Theorem 1]{MartinRosierRouchon-2016} that we recall here:
\begin{theorem}{\cite[Theorem 1]{MartinRosierRouchon-2016}}
	\label{Thm-MartinRosierRouchon}
		Let $u \in \mathscr{R}_L$. Then $u$ can be extended analytically to the set $ \mathcal{S}(L)$. In other words, 
		\begin{equation}
			\label{Supset-Reachable}
			 \mathscr{R}_L \subset \mathscr{A}(L).
		\end{equation}
\end{theorem}
In fact, \cite{MartinRosierRouchon-2016} also describes a subset of $\mathscr{R}$ in terms of the analytic extension of the functions, but shows the following weaker form of \eqref{Subset-Reachable}: If $u \in L^2(-L,L)$ can be extended analytically to the ball $B(0,R)$ for $R > e^{(2e)^{-1}} L$ ( $\simeq 1.2 L$), then $u \in \mathscr{R}_L$. 
\\
The article \cite{MartinRosierRouchon-2016} is to our knowledge the only one describing the reachable set $\mathscr{R}$ without the use of the eigenfunctions of the Laplace operator. If one uses the basis of eigenfunctions of the Laplace operator, simply given by $(\sin(n \pi x/L))_{n \geq 1}$, the results of \cite{ErvZuazuaARMA} (which is a slightly more precise version of \cite{FattoriniRussel71} in this 1d case) yield:
\begin{equation}
	\label{Subset-ErvZua}
	\left\{
		u(x) = \sum_{n \geq 1} c_n \sin\left(\frac{n \pi (x+L)}{2L} \right)  \ 
		\hbox{ such that }
		 \ \sum_n |c_n|^2 n e^{  n \pi }< \infty
	\right\}	
	\subset \mathscr{R}_L.
\end{equation}
Note that condition $	\sum_n |c_n|^2 n e^{n \pi }< \infty$ implies that the function $u(x) = \sum_{n \geq 1} c_n \sin(n \pi (x+L)/2L)$ admits an analytic extension in the strip $\{(x + \ic y)\ |\ |y | < L \}$. Besides, it also implies the boundary conditions $u(-L) = u (L) = 0$ and for all $ n \geq 1$, $ (\partial_{xx})^n u(-L) = (\partial_{xx})^n u(L) = 0$. As pointed out in \cite{MartinRosierRouchon-2016}, this latter condition is rather conservative and should not be relevant as the control is acting on the boundary $x\in \{-L, L\}$.
\\
The proof of Theorem \ref{Thm-Main} will be presented in Section \ref{Sec-Proof-Main}, and is inspired by a Carleman type inequality for the adjoint equation (after the change of variable $t \to T - t$ and having done the formal limit $T \to \infty$):
\begin{equation}
	\label{Eq-Heat-1d-z}
		\left\{
			\begin{array}{ll}
				\partial_t z - \partial_{xx} z = 0 & \hbox{ in } (0,\infty) \times (-L,L), 
				\\
				z(t,-L) = z(t,L) = 0  & \hbox{ in } (0,\infty), 
				\\
				z(0, x) = z_0(x) & \hbox{ in } (-L,L).
			\end{array}
		\right.
\end{equation}
Namely, let us recall that the work \cite{ErvZuazuaARMA} proves the following observability type estimate: There exists $C>0$ such that any smooth solution $z$ of \eqref{Eq-Heat-1d-z} satisfies:
\begin{equation}
	\label{Obs-1-d-direct-Old}
		\int_0^\infty \int_{-L}^L |z(t,x)|^2 \exp\left(- \frac{L^2}{2t} \right)\, dt \, dx \leq C \int_0^\infty  \left(|\partial_x z(t,-L)|^2 + |\partial_x z(t,L)|^2 \right)\, dt.
\end{equation}
In Section \ref{Sec-Proof-Carleman}, we will prove the following improved version of \eqref{Obs-1-d-direct-Old}:
\begin{theorem}
	\label{Thm-Carleman}
	For all $T >0$ satisfying
	\begin{equation}
		\label{Time-Condition-L}
		\pi T > L^2, 
	\end{equation}
	there exists a constant $C>0$ such that for any smooth solution $z$ of \eqref{Eq-Heat-1d-z}, we have the observability inequality:
	\begin{multline}
		\label{Obs-1-d-direct}
		\int_{-L}^L |z(T,x)|^2 \exp\left(\frac{x^2 - L^2}{2T} \right)\, dx 
		+ 
		\int_0^\infty \int_{-L}^L |z(t,x)|^2 \exp\left(\frac{x^2 - L^2}{2t} \right)\, dt \, dx 
		\\
		\leq C \int_0^T t \left( |\partial_x z(t, -L)|^2 +  |\partial_x z(t, L)|^2 \right) \, dt.
	\end{multline}
\end{theorem}
We emphasize that the improvement of \eqref{Obs-1-d-direct} with respect to \eqref{Obs-1-d-direct-Old} is due to the presence of the weight function depending on $x$ in \eqref{Obs-1-d-direct}. Besides, as we will see in Section \ref{Sec-Proof-Carleman}, the proof of Theorem \ref{Thm-Carleman} is more direct than the proof of \eqref{Obs-1-d-direct-Old} in \cite{ErvZuazuaARMA} as it is not based on the observability of the corresponding wave operator.
\\
 In fact, our proof of \eqref{Obs-1-d-direct} closely follows the one of the classical Carleman estimates for the heat equation derived for instance in \cite{FursikovImanuvilov}. In that context, the corresponding weight function $\exp(( x^2- L^2)/4t)$ corresponds to the inverse of the exponential envelop of the kernel 
 $$
 	k_{L}(t,x) = \frac{1}{\sqrt{4 \pi t}} \sin\left( \frac{x L}{2t} \right) \exp\left( \frac{L^2 - x^2}{4t} \right), 
$$
which corresponds to a solution of \eqref{Eq-Heat-1d-z}$_{(1)}$ (in fact, it is the usual Gaussian but translated in the complex plane $x \mapsto x + i L$) used in the transmutation technique in \cite{ErvZuazuaARMA}. Furthermore, this function $k_L$ can be used to check that estimate \eqref{Obs-1-d-direct} is sharp with respect to the blow up of the weight close to $t= 0$.
\\
The condition \eqref{Time-Condition-L} appears naturally in our proof of \eqref{Obs-1-d-direct}. One could naturally think that this condition is remanent from some kind of parabolic version of Ingham's inequality (\cite{Ing}). But this is not the case.  In fact, condition \eqref{Time-Condition-L} rather comes from the fact that, when applying \eqref{Obs-1-d-direct-Old} to the solution $\exp(- \pi^2 t/(4L^2)) \sin( \pi (x+L)/2L)$ of \eqref{Eq-Heat-1d-z}, the weight function in time appearing is $\exp(- 2 \pi^2 t/ (4 L^2) - L^2/2t )$, whose monotony changes precisely at $T^* = L ^2/ \pi$.
\\�\par
One then needs to interpret Theorem \ref{Thm-Carleman} in terms of a dual controllability statement. This mainly consists in the usual duality statement between controllability and observability of the adjoint equation (see e.g. \cite{DoleckiRussell,LionsSIAM88}). To be more precise, we obtain the following result (see Subsection \ref{Subsec-Lem-Duality}):
\begin{lemma}
	\label{Lem-Duality}
	Let $g \in L^2(0, \infty; L^2(-L,L))$ be such that 
	\begin{equation}
		\label{Assumption-g-1}
		\int_0^{\infty} \int_{-L}^L |g(t,x)|^2 \exp\left( \frac{L^2 - x^2}{2t} \right) \, dx \, dt
	< \infty,
	\end{equation}
	and $T$ satisfying \eqref{Time-Condition-L}. If $w$ denotes the solution of 
	\begin{equation}
	\label{Eq-Heat-1d-w-0-L}
		\left\{
			\begin{array}{ll}
				- \partial_t w - \partial_{xx} w = g & \hbox{ in } (0,T) \times (-L,L), 
				\\
				w(t,-L) = w(t,L) = 0  & \hbox{ in } (0,T), 
				\\
				w(T,x) = w_T (x) & \hbox{ in } (-L,L),
			\end{array}
		\right.
	\end{equation}
	for some $w_T \in L^2(-L,L )$, 
	then $w_0(x) = w(0,x)$ belongs to the reachable set $\mathscr{R}_L$.
\end{lemma}
However, it is not completely straightforward to use Lemma \ref{Lem-Duality} as the fundamental solution of the heat equation in a bounded domain involves a discrete summation (namely, the method of images yields a fundamental solution under the form of a sum of odd and even translations of the usual Gaussian kernel). Instead, we prefer to rely on the following result, proved in Section \ref{Sec-Duality}: 
\begin{theorem}
	\label{Thm-Reachable-In-Terms-Of-Source-Term}
	Let ${L_0} >L$ and $\alpha:[-{L_0},{L_0}]\to \mathbb{R}$ be a continuous function on $[-{L_0},{L_0}]$. For $h \in L^2(-{L_0}, {L_0})$, we define 
	\begin{equation}
		\label{Def-g}
		g(t,x) = \frac{1}{t^{3/2}} \exp\left( \frac{x^2 - {L_0}^2}{4t} + \ic  \frac{\alpha(x)}{4t} \right) h(x), 
		\quad (t,x) \in (0,\infty) \times (-{L_0}, {L_0}). 
	\end{equation}
	Then the state $w_0$ defined on $(-L, L)$ by 
	\begin{equation}
		\label{Def-w-0}
		w_0(x) 
		= 
		\int_0^\infty \int_{-{L_0}}^{{L_0}} \frac{1}{ \sqrt{4\pi t}} \exp\left( - \frac{(x - \tilde x)^2}{4t} \right) g(t,\tilde x) \, dt d\tilde x, 
		\quad 
		x \in (-L, L), 
	\end{equation}
	is well-defined and belongs to the reachable set $\mathscr{R}_L$.
	\\
	Besides, $w_0$ can alternatively be written as
	\begin{equation}
		\label{Def-w-0-2}
		w_0(x) =   \frac{2}{\sqrt{\pi}} \int_{-{L_0}}^{{L_0}} \frac{h(\tilde x)}{(x - \tilde x)^2 + {L_0}^2 - {\tilde x}^2 + \ic \alpha(\tilde x)}  \, d\tilde x, 
	\quad x \in (-L,L). 
	\end{equation}
\end{theorem}	
Under the conditions of Theorem \ref{Thm-Reachable-In-Terms-Of-Source-Term}, explicit computations yield that for any $\alpha \in C^0([-{L_0},{L_0}]; \mathbb{R})$, the range of the operator $K_\alpha: L^2(-{L_0},{L_0}) \to L^2(-L,L)$ given for $h \in L^2(-{L_0}, {L_0})$ by 
\begin{equation}
	\label{Def-K-alpha}
	K_\alpha (h)(x) =   \frac{2}{\sqrt{\pi}} \int_{-{L_0}}^{{L_0}} \frac{h(\tilde x)}{(x - \tilde x)^2 + {L_0}^2 - {\tilde x}^2 + \ic \alpha(\tilde x)}  \, d\tilde x, 
	\quad x \in (-L,L), 
\end{equation}
is contained in the reachable set $\mathscr{R}_L$. Therefore, Theorem \ref{Thm-Reachable-In-Terms-Of-Source-Term} can be rewritten as 
\begin{equation}
	\label{Range-K-Alpha}
	\forall \alpha \hbox{ satisfying the assumptions of Theorem \ref{Thm-Reachable-In-Terms-Of-Source-Term}}, 
	\quad \hbox{Range} (K_\alpha) \subset \mathscr{R}. 
\end{equation}
The proof of Theorem \ref{Thm-Main} will thus mainly reduced to choosing carefully the above functions $\alpha$ and $h$ in Theorem \ref{Thm-Reachable-In-Terms-Of-Source-Term} so that all function which admits an analytic extension in a set $\mathcal{S}(L_0)$ for $L_0 > L$ can be decomposed into a finite sum of elements in the images of the above operators $K_\alpha$. This property will be achieved by using Cauchy's formula for holomorphic functions among contours which coincide with the singularities of the kernel of the above operator $K_\alpha$. Details of the proof are given in Section \ref{Sec-Proof-Main}.
\\
As it turns out, see Section \ref{Sec-Proof-Main}, we will require the use of non-trivial functions $\alpha$ in a critical way. This might be surprising at first as this function introduces strong time oscillations in the source term of the heat equation in \eqref{Def-g}. In other words, we need these strong oscillations to reach the whole reachable set. 
\\ \par
\noindent{\bf Scientific Context.} The characterization of the reachable set of the heat equation is a rather old issue, whose study probably started with the pioneering work \cite{FattoriniRussel71} studying this question in dimension one using harmonic analysis techniques. The result of \cite{FattoriniRussel71} was then slightly improved into \eqref{Subset-ErvZua} in \cite{ErvZuazuaARMA} using the so-called {\it transmutation} technique allowing to write solutions of the wave equations in terms of solutions of the heat equation \eqref{Eq-Heat-1d-z}.
\\
More recently, P. Martin, L. Rosier and P. Rouchon proposed in the work \cite{MartinRosierRouchon-2016} to characterize the reachable set of the heat equation in the 1d case by a description on the set on which the reachable states are analytic. As explained above, this description yields that if a state admits an analytic extension on the ball $B(0,R)$ for $R > e^{(2e)^{-1}} L$ ($\simeq 1.2 L$), then it belongs to the reachable set $\mathscr{R}_L$. The approach in \cite{MartinRosierRouchon-2016} relies on the {\it flatness approach}, which has been developed recently by P. Martin, L. Rosier and P. Rouchon, see \cite{MartinRosierRouchon-Automatica-2014,MartinRosierRouchon-Sicon-2016}.
\\
Of course, describing the reachable set of the heat equation is also related to the results of controllability to trajectories for the heat equation, which by linearity are equivalent to the results on null-controllability. In this context, the breakthrough came from the introduction of Carleman estimates to obtain observability results for the heat equation in any dimension from basically any non-open subset, see \cite{FursikovImanuvilov,LebRob}. Nonetheless, in general, Carleman estimates are not suitable to provide sharp estimates as one has very little control on the coefficients appearing in them. Theorem \ref{Thm-Carleman} is a very specific case in which the parameters can be explicitly computed. 
\\
With that in mind, one could also relate the Carleman estimate in Theorem \ref{Thm-Carleman} with the Hardy Uncertainty principle obtained in \cite{Escauriaza-Kenig-Ponce-Vega-2016}. In some sense, the weight that we are using is a {\it limiting Carleman weight}, in the sense that the conjugated operator appearing in the proof of Theorem \ref{Thm-Carleman} satisfies a degenerate convexity condition, see Remark \ref{Rem-Degenerate}.
\\
Let us also emphasize that there are several works related to the cost of controllability of the heat equation in short time. Let us quote in particular the works by \cite{Miller04a,Miller06a,Miller06b,TenenbaumTucsnak07} studying these questions. It was thought for a while that the understanding of the blow up of the controllability of the heat equation in short time would be more or less equivalent to a good characterization of the reachable set, but this was recently disproved in \cite{Lissy-2015}. We refer the interested reader to this latter work for a more detailed discussion on this fact. 
\\
These issues are also related to the questions raised in \cite{CoronGuerrero} concerning the controllability of a viscous transport equation with vanishing viscosity parameter, see \cite{LissyCras2012,Lissy-ErvZua-Cor-Gue}. In that sense, our work suggests that observability results for the heat equation stated only in terms of their spectral decomposition could possibly be reinforced by considering space weighted functional settings appropriate to the control problem at hand.
\\ \par
\noindent{\bf Outline.} This article is organized as follows. Section \ref{Sec-Proof-Carleman} gives the proof of Theorem \ref{Thm-Carleman}. Section \ref{Sec-Duality} establishes Theorem \ref{Thm-Reachable-In-Terms-Of-Source-Term}. In Section \ref{Sec-Proof-Main}, we prove Theorem \ref{Thm-Main}. We finally provide the reader with further comments in Section \ref{Sec-Further}.
\\ \par
\noindent{\bf Acknowledgements.} The authors deeply thank Michel Duprez and Pierre Lissy for several useful comments related to this work.
\section{Proof of Theorem \ref{Thm-Carleman}}\label{Sec-Proof-Carleman}
	Let $z$ be a smooth solution of \eqref{Eq-Heat-1d-z}, and introduce the new unknown (the conjugated variable):
	\begin{equation}
		\label{tilde-z-z}
		\tilde z(t,x) = z(t,x) t \exp\left( \frac{x^2 - L^2}{4t} \right), \quad (t,x ) \in (0, \infty) \times (-L, L). 
	\end{equation}
	It satisfies the equations
	\begin{equation}
		\label{Eq-Conjugated-Heat-Free-1d}
		\left\{
			\begin{array}{ll}
				\ds \partial_t \tilde z + \frac{x}{t} \partial_x \tilde z - \frac{1}{2t} \tilde z- \partial_{xx} \tilde z - \frac{L^2}{4t^2} \tilde z = 0, \quad & (t,x) \in (0, \infty) \times (-L,L), 
				\smallskip\\
				\tilde z(t,-L) = \tilde z(t,L) = 0, \quad & t \in (0, \infty), 
				\\
				\tilde z(0, x) = 0, \quad   & x \in( -L,L).
			\end{array}
		\right.
	\end{equation}
	We then introduce the energy $E(t)$ and the dissipation $D(t)$ defined for $t>0$ by 
	\begin{eqnarray}
		\label{Def-E}
		E(t) &=& \int_{-L}^L | \tilde z(t,x)|^2 \, dx,
		\\
		\label{Def-D}
		D(t) &=& \int_{-L}^L |\partial_x \tilde z(t,x)|^2 \, dx - \frac{L^2}{4t^2} \int_{-L}^L |\tilde z(t,x)|^2 \, dx. 
	\end{eqnarray}
	Easy computations show that they satisfy the following ODEs: for all $t >0$, 
	\begin{align}
		\label{Evol-E}
		&
		\frac{dE}{dt}(t) - \frac{2}{t} E(t) + 2 D(t) = 0, 
		\\
		\label{Evol-D}
		&
		\frac{dD }{dt} (t) + 2 \int_{-L}^L \left| - \partial_{xx} \tilde z(t,x) - \frac{L^2}{4t^2} \tilde z(t,x) \right|^2 \, dx = \frac{L}{t} \left( |\partial_x \tilde z(t,-L)|^2 +|\partial_x \tilde z(t,L)|^2  \right).
	\end{align}
	But, by Poincar\'e's inequality, $D(t)$ is non-negative for $t\geq L^2/\pi$. Let then $T> L^2/\pi$ as in \eqref{Time-Condition-L}. Integrating \eqref{Evol-D} between $0$ and $T$ and using the fact that $D(0) = 0$ due to vanishing behavior of the weight function close to $t = 0$, we get 
	\begin{equation}
		\label{Intermediate-Step-Carleman}
		\int_0^T \int_{-L}^L \left| - \partial_{xx} \tilde z(t,x) - \frac{L^2}{4t^2} \tilde z(t,x) \right|^2 \, dtdx 
		\leq \int_0^T \frac{L}{2t} \left( |\partial_x \tilde z(t,-L)|^2 +|\partial_x \tilde z(t,L)|^2  \right) \, dt.
	\end{equation}
	But $\tilde z$ satisfies the boundary conditions $\tilde z(t,-L) = \tilde z(t,L) = 0$ for all $t >0$. Therefore, multiplying $- \partial_{xx} \tilde z- \frac{L^2}{4t^2} \tilde z $ by $ 2 x \overline{\partial_x \tilde z}$, for all $ t >0$  we obtain 
	\begin{align*} 
		& \int_{-L}^L |\partial_x \tilde z(t,x)|^2 \, dx + \frac{L^2}{4t^2} \int_{-L}^L |\tilde z(t,x)|^2 \, dx 
		\\
		&
		= 
		2 \Re\left( \int_{-L}^L \left( - \partial_{xx} \tilde z(t,x) - \frac{L^2}{4t^2} \tilde z(t,x) \right) x \overline{\partial_x \tilde z} \, dx \right)
		+ L  \left(  |\partial_x \tilde z(t,-L)|^2  + |\partial_x \tilde z(t,L)|^2 \right)
		\\
		& 
		\leq
		L^2  \int_{-L}^L \left| - \partial_{xx} \tilde z(t,x) - \frac{L^2}{4t^2} \tilde z(t,x) \right|^2 \, dx
		+ 
		\int_{-L}^L |\partial_x \tilde z(t,x)|^2 \, dx 
		+
		 L  \left(  |\partial_x \tilde z(t,-L)|^2  + |\partial_x \tilde z(t,L)|^2 \right), 
	\end{align*}
	so that for all $t >0$, 
	$$
		\frac{L^2}{4t^2} \int_{-L}^L |\tilde z(t,x)|^2 \, dx 
		\leq 
		L^2  \int_{-L}^L \left| - \partial_{xx} \tilde z(t,x) - \frac{L^2}{4t^2} \tilde z(t,x) \right|^2 \, dx
		+
		 L  \left(  |\partial_x \tilde z(t,-L)|^2  + |\partial_x \tilde z(t,L)|^2 \right). 
	$$
	Using this last estimate in \eqref{Intermediate-Step-Carleman}, we derive
	\begin{equation}
		\label{Est-on-0-T}	
		\int_0^T \int_{-L}^L \frac{1}{t^2} |\tilde z(t,x)|^2 \, dt dx \leq C \int_0^T \frac{1}{t} \left(  |\partial_x \tilde z(t,-L)|^2  + |\partial_x \tilde z(t,L)|^2 \right)\, dt,  
	\end{equation}
	Besides, from \eqref{Evol-E}, again using Poincar\'e estimate, for all $t \geq  T > L^2/ \pi$, 
	\begin{equation}
		\label{ODE-on-T-infinity}	
		\frac{d}{dt} \left( \frac{E(t)}{t^2} \right) + \frac{L^2}{T^2} \left(\frac{\pi^2 T^2}{L^4} -1 \right) \frac{E(t)}{t^2} \leq 0, 
	\end{equation}
	while $t \mapsto E(t)/t^2$ is decreasing on $(L^2/\pi, T)$ from \eqref{Evol-E}:
	\begin{equation}
		\label{Est-at-T}	
		\frac{E(T)}{T^2} 
		\leq 
		\frac{1}{ T - L^2/\pi } \int_{L^2/\pi}^T \frac{E(t)}{t^2} \, dt
		\leq
		\frac{1}{ T - L^2/\pi } \int_0^T  \int_{-L}^L \frac{1}{t^2} |\tilde z(t,x)|^2 \, dt dx.
	\end{equation}
	Therefore, combining \eqref{Est-on-0-T}--\eqref{ODE-on-T-infinity}-\eqref{Est-at-T}, we easily derive
	$$
		\frac{E(T)}{T^2} + \int_0^\infty \int_{-L}^L \frac{1}{t^2} |\tilde z(t,x)|^2 \, dt dx \leq C \int_0^T \frac{1}{t} \left(  |\partial_x \tilde z(t,-L)|^2  + |\partial_x \tilde z(t,L)|^2 \right)\, dt,  
	$$
	Using \eqref{tilde-z-z}, we immediately obtain \eqref{Obs-1-d-direct}.
\begin{remark}
	\label{Rem-Degenerate}
	Note that the conjugated operator in \eqref{Eq-Conjugated-Heat-Free-1d} is in some sense degenerate. Indeed, the conjugated operator in \eqref{Eq-Conjugated-Heat-Free-1d} is 
	$$
		\partial_t + \frac{x}{t} \partial_x - \frac{1}{2t} - \partial_{xx} - \frac{L^2}{4t^2}, 
	$$	
	which can be written as $A+B+R$ with 
	$$
		A = \partial_t + \frac{x}{t} \partial_x + \frac{1}{2t}, 
		\quad 
		B = - \partial_{xx} - \frac{L^2}{4t^2}, 
		\quad \hbox{�and }�
		R = - \frac{1}{t},
	$$
	and which satisfies
	$$
		A^* = - A, \quad B^* = B, \quad [A, B] = - \frac{2}{t} B, 
	$$
	while the operator $R$ is of lower order.
	\\
	In particular, if the symbol of the operator $A$ and of the operator $B$ cancels, the symbol of the commutator $[A,B]$ vanishes as well. The convexity condition needed to get Carleman estimates is therefore degenerate with our choice of weights. Such degenerate weights have appeared in the literature lately in the context of the Calder\'on problem, see in particular the works on the limiting Carleman weights, see \cite{KSUhlmann-Calderon-partialdata} and subsequent works.
\end{remark}
\section{Duality results}\label{Sec-Duality}
\subsection{Proof of Lemma \ref{Lem-Duality}}\label{Subsec-Lem-Duality}
	Let $w$ denote the solution of \eqref{Eq-Heat-1d-w-0-L} with source term $g$ satisfying \eqref{Assumption-g-1} and initial datum $w_T \in L^2(-L,L)$. 
	\\
	For the proof below and in order to simplify the notations, we further assume that $g$ and $w_T$ are real-valued. This can be done without loss of generality  by applying the result to the real part of $(g, w_T)$ and the imaginary part of $(g,w_T)$ in case $(g,w_T)$ are complex-valued.
	\\
	Let us then define the functional $J$ as follows: for $z_0 \in L^2(-L,L; \mathbb{R})$, 
	\begin{equation}
		\label{Def-J-z-0}
		{J}(z_0) 
		:= 
		\frac{1}{2} \int_0^T  \left( \vert \partial_x z(t,-L) \vert^2 + 	\vert \partial_x z(t,L) \vert^2 \right) \, dt 
		- \int_0^{+\infty} \int_{-L}^L g \, z \, dt\, dx
		- \int_{-L}^L z(T,x) w_T(x) \,dx, 
	\end{equation}
	where $z$ denotes the solution of \eqref{Eq-Heat-1d-z} with initial datum $z_0$.
	\\
	From the Carleman estimate in Theorem \ref{Thm-Carleman} and assumption \eqref{Assumption-g-1}, the functional $J$ can be extended by continuity on the space 
	\begin{equation}
		X = \overline{\{ z_0 \in L^2(-L,L; \mathbb{R})\}}^{\norm{\cdot}_{obs}}, 
		\hbox{ with } 
		\norm{z_0}_{obs}^2 
		= 
		\int_0^T  \left( \vert \partial_x z(t,-L) \vert^2 + 	\vert \partial_x z(t,L) \vert^2 \right) \, dt, 
	\end{equation}
	where $z$ denotes the corresponding solution of \eqref{Eq-Heat-1d} with initial datum $z_0$. Note that using \eqref{Obs-1-d-direct}, one can associate to $z_0 \in X$ a solution $z$ of \eqref{Eq-Heat-1d-z}$_{(1,2)}$ and normal traces $\partial_x z(t, \pm L)$ with $z \exp((x^2 - L^2)/4t) \in L^2(0,\infty; L^2(-L,L))$ and $ \partial_x z(t,\pm L) \in L^2(0, T)$.
	\\
	Besides, the functional $J$ is strictly convex and coercive on $X$ from \eqref{Obs-1-d-direct}. Therefore, it admits a unique minimizer $Z_0$ in $X$. Let us denote by $\partial_x Z(t,-L),\, \partial_x Z(t,L)$ the corresponding normal traces, and set 
	\begin{equation}
		v_- (t) =  \partial_x Z(t,-L), 
		\quad
		v_+ (t) =  - \partial_x Z(t,L), 
		\quad \hbox{ in } (0,T).
	\end{equation}
	Using then that $J(Z_0) \leq J(0)$ and the Carleman estimate \eqref{Obs-1-d-direct}, one easily checks that 
	\begin{multline}
		\norm{v_-(t)}_{L^2(0,T)}^2 
		+ 
		\norm{v_+(t)}_{L^2(0,T)}^2
		\\
		\leq
		C 
		\int_0^{\infty} \int_{-L}^L |g(t,x)|^2 \exp\left( \frac{L^2 - x^2}{2t} \right) \, dx \, dt
		+ 
		C\int_{-L}^L |w_T(x)|^2 \exp\left( \frac{L^2 - x^2}{2T} \right)  \, dx.
		\label{Norm-v-L2}
	\end{multline}
	Furthermore, the Euler-Lagrange equation of $J$ at $Z_0$ in the direction $z_0 \in L^2(-L,L; \mathbb{R})$ yields:
	\begin{equation}
		\label{Euler-Lag-z-0}
		0 =  \int_0^T \left( v_-(t) \partial_x z(t,-L) - v_+(t)\partial_x z(t,L) \right) \, dt 
		- \int_0^{+\infty} \int_{-L}^L g \, z \, dt\, dx
		- \int_{-L}^L z(T,x) w_T(x) \,dx.
	\end{equation}	
	But, multiplying the equation \eqref{Eq-Heat-1d-w-0-L} satisfied by $w$ by $z$ solution of \eqref{Eq-Heat-1d-z}, we get the identity:
	\begin{equation}
		\label{Duality-Paring-w-z}
		\int_0^{+\infty} \int_{-L}^L g \, z \, dt\, dx
		+ 
		\int_{-L}^L z(T,x) w_T(x) \,dx
		= 
		\int_{-L}^L z_0(x) w(0,x) \,dx.
	\end{equation}
	Therefore identity \eqref{Euler-Lag-z-0} can be rewritten as follows: for all $z_0 \in L^2(-L,L; \mathbb{R})$, denoting by $z$ the solution of \eqref{Eq-Heat-1d-z}, one has
	\begin{equation}
		\label{w-0-v-pm}
		\int_{-L}^L z_0(x) w(0,x) \,dx
		= 
		 \int_0^T \left( v_-(t) \partial_x z(t,-L) - v_+(t)\partial_x z(t,L) \right) \, dt .
	\end{equation}
	Recall then that $v_-, v_+$ belong to $L^2(0,T)$ according to \eqref{Norm-v-L2}, and let us then define $u$ the solution of 
	\begin{equation}
		\label{Eq-u-Dual}
		\left\lbrace
			\begin{array}{ll}
				\partial_t u + \partial_{xx} u = 0 & \text{ in } (0,T) \times (-L,L) ,
				\\
				u(t,- L) =  v_{-} (t)  & \text{ in } (0,T),
				\\
				u(t,- L) =  v_{+} (t)  & \text{ in } (0,T),
				\\
				u(T,x) = 0 & \text{ in } (-L,L).
			\end{array}
		\right.
	\end{equation}
	We claim that $u(0, \cdot) = w(0,\cdot)$. Indeed, for $z_0 \in L^2(-L,L; \mathbb{R})$, if we multiply the equation of $u$ in \eqref{Eq-u-Dual} by the solution $z$ of \eqref{Eq-Heat-1d-z}, we get
	\begin{equation}
		\label{u-0-v-pm}
		\int_{-L}^L u(0, x) z_0(x) \, dx
		= 
		- \int_0^T v_+(t) \partial_x z(t,L)\, dt 
		+\int_0^T v_-(t) \partial_x z(t,-L)\, dt . 
	\end{equation}
	Therefore, comparing \eqref{w-0-v-pm} with \eqref{u-0-v-pm}, we get that 
	$$
		\forall z_0 \in L^2(-L,L; \mathbb{R}), \quad \int_{-L}^L (u(0,x) - w(0,x)) z_0(x) \, dx = 0, 
	$$
	that is $u(0,\cdot) = w(0, \cdot)$. 
	\\
	We then simply remark that doing the change of unknowns $\tilde u(t, \cdot) = u(T-t, \cdot)$, $\tilde v_{\pm}(t) = v_\pm(T-t)$, $u(0, \cdot) = \tilde u(T, \cdot)$ is a reachable state for \eqref{Eq-Heat-1d} with controls $\tilde v_\pm$, i.e. $ u(0, \cdot) \in  \mathscr{R}_L$. As $u(0,\cdot) = w(0, \cdot)$, we have thus obtained that $w(0, \cdot) \in \mathscr{R}_L$.
\subsection{Proof of Theorem \ref{Thm-Reachable-In-Terms-Of-Source-Term}}
	Set $g$ as in \eqref{Def-g} and define, for $(t,x) \in (0,T] \times (-{L_0}, {L_0})$, the function $w(t,x)$ as follows
	\begin{eqnarray}
		\label{Def-w-t-x}
		w(t,x) 
		& = &
		\int_t^\infty \int_{-{L_0}}^{{L_0}} \frac{1}{ \sqrt{4\pi (\tilde t - t)}} \exp\left( - \frac{(x - \tilde x)^2}{4(\tilde t -t)} \right) g(\tilde t,\tilde x) \, d\tilde t d\tilde x
		\\
		\label{Def-w-t-x-2}
		& = & 
		\int_0^\infty \int_{-{L_0}}^{{L_0}} \frac{1}{ \sqrt{4\pi s}} \exp\left( - \frac{(x - \tilde x)^2}{4s} \right) g(s+t,\tilde x) \, ds d\tilde x.			
	\end{eqnarray}
	Our goal is to check that $w$ solves 
	\begin{equation}
	\label{Eq-Heat-1d-w}
		\left\{
			\begin{array}{ll}
				- \partial_t w - \partial_{xx} w = g & \hbox{ in } (0,T) \times (-L,L), 
				\\
				w(t,-L) = v_-(t) & \hbox{ in } (0,T), 
				\\
				w(t,L) = v_+(t)  & \hbox{ in } (0,T), 
				\\
				w(T,x) = w_T (x) & \hbox{ in } (-L,L),
			\end{array}
		\right.
	\end{equation}
	with appropriate choice of functions $w_T \in L^2(-L,L)$, $v_- \in L^2(0,T)$ and $v_+ \in L^2(0,T)$ and that $w_0$ defined in \eqref{Def-w-0} simply is the trace of $w$ at time $t = 0$.
	\\
	Indeed, if \eqref{Eq-Heat-1d-w} holds, we can decompose $w$ as $w = \tilde w + \hat w$, with $\tilde w$ satisfying \eqref{Eq-Heat-1d-w-0-L} with source term $g$ and with initial condition $\tilde w(T) = w_T$ and $\hat w$ satisfying the equation 
	\begin{equation}
	\label{Eq-Heat-1d-w-hat}
		\left\{
			\begin{array}{ll}
				- \partial_t \hat w - \partial_{xx} \hat w = 0 & \hbox{ in } (0,T) \times (-L,L), 
				\\
				\hat w(t,-L) = v_-(t) & \hbox{ in } (0,T), 
				\\
				\hat w(t,L) = v_+(t)  & \hbox{ in } (0,T), 
				\\
				\hat w(T,x) = 0 & \hbox{ in } (-L,L),
			\end{array}
		\right.
	\end{equation}
	for which one immediately has (by the change of variables $t \to T-t$) that $\hat w(0, \cdot) \in \mathscr{R}_L$. The state $\tilde w(0, \cdot)$ belongs to $\mathscr{R}_L$ since $g$ defined in \eqref{Def-g} satisfies \eqref{Assumption-g-1} due to the condition ${L_0} > L$ and so Lemma \ref{Lem-Duality} applies. This eventually implies that $w(0, \cdot)$ belongs to $\mathscr{R}_L$ as $\hat w(0, \cdot)$ and $\tilde w(0, \cdot)$ belong to $\mathscr{R}_L$.
	\\
	We therefore first focus on the proof of the fact that $w$ in \eqref{Def-w-t-x} satisfies \eqref{Eq-Heat-1d-w} with $w_T \in L^2(-L,L)$, $v_- \in L^2(0,T)$ and $v_+ \in L^2(0,T)$ and that $w_0$ defined in \eqref{Def-w-0} simply is the trace of $w$ at time $t = 0$.
	\\ \par
	Let us now prove that $w$ in \eqref{Def-w-t-x} satisfies \eqref{Eq-Heat-1d-w}.
	\\
	We first remark that $g$ in \eqref{Def-g} satisfies, $for all (t,x) \in (0, \infty) \times (-{L}, {L})$
	\begin{equation}
		\label{Est-g}
		|g(\tilde t,\tilde x)| 
		\leq 
		\frac{1}{\tilde t^{3/2}} \exp\left( \frac{\tilde x^2 - {L_0}^2}{4 \tilde t}\right) |h(\tilde x)| 
		\leq  \frac{1}{\tilde t^{3/2}} |h(\tilde x)| , 
	\end{equation}
	with $h \in L^2(-{L_0}, {L_0})$. The continuity of $w$ in \eqref{Def-w-t-x-2} is therefore easy to prove on all sets of the form $(t,x) \in (\varepsilon, \infty) \times (-{L_0}, {L_0})$ with $\varepsilon >0$, as the decay in $\tilde t$ in \eqref{Est-g} makes the integral convergent for $s= \tilde t -t $ close to infinity while the integrability for $s$ close to $0$ simply comes from the integrability of $s \mapsto s^{-1/2}$ close to $0$. But the continuity close to $t= 0$ is more delicate to obtain. We will simply show that $w$ in \eqref{Def-w-t-x-2} is continuous on $(t,x) \in (0, \infty) \times (-\ld, \ld)$ for $\ld \in (L, {L_0})$. Indeed, let us set $\ld \in (L, {L_0})$, and let us rewrite $w$ in \eqref{Def-w-t-x-2} as:
	\begin{equation}
		\label{Def-w-t-x-3}
		w(t,x) = \int_0^\infty \int_{-{L_0}}^{{L_0}} \frac{1}{ \sqrt{4\pi s}} \frac{1}{(t+s)^{3/2}} \exp\left( - \frac{(x - \tilde x)^2}{4s} + \frac{\tilde x^2- {L_0}^2}{4 (t+s)} + \ic \frac{\alpha(\tilde x)}{4(t+s)}\right) h(\tilde x) \, ds d\tilde x.
	\end{equation}
	Under this form, it is clear that what matters it the sign of 
	$$
		P(t,s,x,\tilde x) = - \frac{(x - \tilde x)^2}{4s} + \frac{\tilde x^2- {L_0}^2}{4 (t+s)}.
	$$
	But for $t \geq 0$, $s \in (0, \infty)$, $x \in [-\ld,\ld]$ and $\tilde x \in [-{L_0}, {L_0}]$, we have
	\begin{eqnarray*}
		P(t,s,x,\tilde x) 
		& =& 
		\frac{1}{4s (t+s)} \left( - (x - \tilde x)^2 t + s (- (x- \tilde x)^2 + \tilde x^2 - {L_0}^2) \right)
		\\
		& \leq &
		\frac{1}{4 (t+s)} (- x^2 + 2 x \tilde x - {L_0}^2)
		\\
		& \leq &
		\frac{1}{4(t+s)} (- x^2 + 2 |x| {L_0} - {L_0}^2)
		\\
		& \leq &
		- \frac{1}{4(t+s)} (\ld- {L_0})^2.
	\end{eqnarray*}
	One then easily deduces that $w$ in \eqref{Def-w-t-x-3} is continuous on $[0,\infty) \times (-\ld,\ld)$ and that its value at $t = 0$ coincides with the formula \eqref{Def-w-0}.
	\\
	We can then set
	$$
		v_-(t) = w(t,-L)\, \hbox{�in } (0,T), \qquad
		v_+(t) = w(t,L) \, \hbox{�in } (0,T), \qquad
		w_T(x) = w(T,x) \hbox{ in } (-L,L), 
	$$ 	
	for which the previous analysis implies $v_- \in L^2(0,T)$, $v_+ \in L^2(0,T)$ and $w_T \in L^2(-L,L)$. 
	\\
	Finally, the fact that $w$ solves the first equation in \eqref{Eq-Heat-1d-w} obviously comes from the fact that the kernel appearing in \eqref{Def-w-t-x} is the heat kernel. 
	\\ \par
	In order to prove formula \eqref{Def-w-0-2}, we simply use Fubini's theorem:
	\begin{eqnarray*}
		w_0(x) 
		& = & 
		\int_0^\infty \int_{-{L_0}}^{{L_0}} \frac{1}{ \sqrt{4\pi } t^{2}�} \exp\left( - \frac{(x - \tilde x)^2}{4t} + \frac{\tilde x^2 - {L_0}^2}{4t} + \ic \frac{\alpha(\tilde x)}{4t}\right) h(\tilde x) \, dt d\tilde x 
		\\
		& = &
		\frac{1}{\sqrt{4\pi}} \int_{-{L_0}}^{L_0} \left(\int_0^\infty \frac{1}{t^2} \exp\left( - \frac{(x - \tilde x)^2}{4t} + \frac{\tilde x^2 - {L_0}^2}{4t} + \ic \frac{\alpha(\tilde x)}{4t}\right) \, dt \right) h(\tilde x) d\tilde x.
		\\
		 &=&   \frac{2}{\sqrt{\pi}} \int_{-{L_0}}^{{L_0}} \frac{h(\tilde x)}{(x - \tilde x)^2 + {L_0}^2 - {\tilde x}^2 + \ic \alpha(\tilde x)}  \, d\tilde x. 
	\end{eqnarray*}
	This concludes the proof of Theorem \ref{Thm-Reachable-In-Terms-Of-Source-Term}.
\section{Proof of Theorem \ref{Thm-Main}}\label{Sec-Proof-Main}
\subsection{Strategy}
As explained in the introduction, our main objective is to study the range of the operators $K_\alpha$ introduced in \eqref{Def-K-alpha} for good choices of functions $\alpha$. 
\\ \par
To start with, we will focus on the case $\alpha = 0$, and in Section \ref{Subsec-Proof-K-alpha-0} we will prove the following:
\begin{proposition}
	\label{Prop-K-alpha-0}
	Let $L>0$ and ${L_0} >L$, and define the operator $K_{0,{L_0}}: L^2(- {L_0}, {L_0}) \to L^2(-L, L)$ by
	\begin{equation}
		\label{Def-K-0}
		K_{0,{L_0}}(h) (x) = \frac{2}{\sqrt{\pi}} \int_{-{L_0}}^{L_0} \frac{h(\tilde x)}{(x - \tilde x)^2 + {L_0}^2 - \tilde x^2} d \tilde x, \quad x \in (-L, L).
	\end{equation}
	 Then any function $k$ defined on $(-L, L)$ which can be extended analytically on the closure of the ball of size ${L_0}$ belongs to the range of the operator $K_{0,{L_0}}$.
\end{proposition}
Proposition \ref{Prop-K-alpha-0} is proved using the Chebyshev polynomials $(U_n)_{n \in \mathbb{N}}$ of the second kind, that is the sequence of polynomials such that for all $n \in \mathbb{N}$: 
\begin{equation}
	\label{Cheb}
	U_n(\cos \theta) = \frac{\sin((n+1) \theta)}{\sin(\theta)}, \quad \theta \in (-\pi, \pi).
\end{equation}
Indeed, they appear naturally as the generating function for the polynomials $U_n$ is given as follows:
\begin{equation}
	\label{Generatrice-Series}
	\sum_{n\geq 0} x^n U_n(\tilde x) = \frac{1}{1 - 2 x \tilde x + x^2} = \frac{1}{(x-\tilde x)^2 + 1 - \tilde x^2}, \quad x, \, \tilde x \in (-1,1).
\end{equation}
Let us also point out that for $h \in L^2(-{L_0}, {L_0})$, we automatically have that  $K_{0, {L_0}}(h)$ admits an analytic extension on the ball of radius ${L_0}$ by the obvious formula:
\begin{equation}
	\label{K-0-ell-z}
	K_{0,{L_0}}(h) (z) = \frac{2}{\sqrt{\pi}} \int_{-{L_0}}^{L_0} \frac{h(\tilde x)}{(z - \tilde x)^2 + {L_0}^2 - \tilde x^2} d \tilde x, \quad z \in B(0, {L_0}).
\end{equation}
In other words, Proposition \ref{Prop-K-alpha-0} proves that the range of the operator $K_{0, {L_0}}$ is very close of being exactly the functions which can be extended analytically to the ball of radius ${L_0}$. This is in fact rather expected due to the similarity of formula \eqref{K-0-ell-z} with the Poisson kernel appearing when solving the Laplace equation in the ball. Note that Proposition \ref{Prop-K-alpha-0} already improves the result of \cite{MartinRosierRouchon-2016}, which proved that functions which can be extended analytically to balls of radius $e^{1/(2e)} L \simeq 1.2 L$ belong to the reachable set $\mathscr{R}$.
\\ \par
\begin{figure}[ht]
\begin{center}
\includegraphics[width = 0.5\textwidth, height = 0.3\textheight]{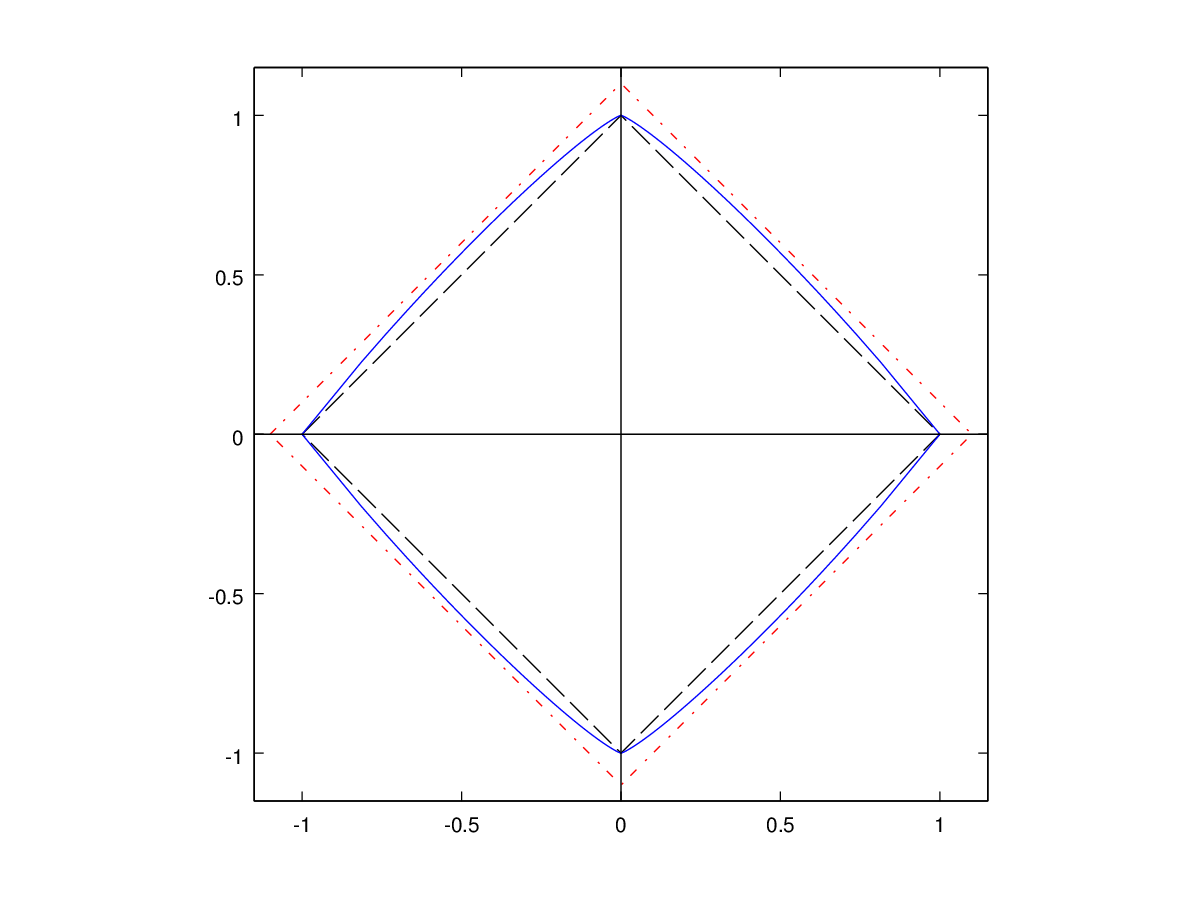}
\end{center}
\caption{In black, $--$ the square $\mathcal{S}(L_0)$ with $L_0 = 1$; in red, $.-$ the square $\mathcal{S}(L_0(1+\varepsilon))$ with $L_0= 1$, $\varepsilon = 0.15$; in blue, a corresponding contour $\mathcal{C}$ given by Lemma \ref{Lem-Choice-Alpha}.}
\label{Fig3}
\end{figure}

The next step then consists in showing that choosing the function $\alpha$ carefully, we can reach any function which can be extended analytically in the neighborhood of the square $\overline{S(L)}$. The basic idea in order to choose the function $\alpha$ appropriately is to remark that the operator $K_\alpha$ in \eqref{Def-K-alpha} has a kernel given for $(x,\tilde x) \in (-L,L) \times (L_0, L_0)$ by 
\begin{multline}
	  \frac{1}{(x - \tilde x)^2 + L_0^2 - {\tilde x}^2 + \ic \alpha(\tilde x)}
	= 
	\frac{1}{X_{\alpha,+}(\tilde x) -X_{\alpha,-}(\tilde x) } \left(\frac{1}{x- X_{\alpha,+}(\tilde x)} 
	- 
	 \frac{1}{x- X_{\alpha,-}(\tilde x)}\right)
	\\
	\hbox{ with } 
	\left\{
		\begin{array}{l}
		\ds X_{\alpha,+}(\tilde x) = \tilde x + \ic \sqrt{L_0^2 - \tilde x^2 + \ic \alpha (\tilde x)},
		\smallskip\\
		\ds  X_{\alpha,-}(\tilde x) = \tilde x - \ic \sqrt{L_0^2 - \tilde x^2 + \ic \alpha (\tilde x)},
		\end{array}	
	\right.
\end{multline}
where we used the complex square-root function cut on the axis $\mathbb{R}_-$.
\\
We claim that it is possible to choose $\alpha$ as follows: 
\begin{lemma}
	\label{Lem-Choice-Alpha}
	Let $L_0 >0$ and $\varepsilon >0$. Then there exists a continuous function $\alpha:(-L_0, L_0) \to \mathbb{R}$ such that:
	\begin{itemize}
		\item[(i)] $\alpha$ is piecewise $C^1(-L_0, L_0)$ and $\alpha$ can be extended as a $C^1$ function on the interval $[-L_0, 0]$ and $[0,L_0]$.
		\item[(ii)] For all $\tilde x \in [-L_0, L_0]$, $X_{\alpha,+}(-\tilde x) = - X_{\alpha,-}(\tilde x)$. 
		\item[(iii)] The set $\{X_{\alpha,-}(\tilde x), \, \tilde x \in [0,L_0]\}$ describes a path included in the set $\{x + \ic y, \, x \geq 0, \, y \leq 0\}\setminus B(0, L_0 )$.
		\item[(iv)] The set $\{X_{\alpha,+}(\tilde x), \, \tilde x \in [0,L_0]\}$ describes a path included in the set $\mathcal{S}(L_0(1+\varepsilon)) \setminus \mathcal{S}(L_0)\cap \{x + \ic y, \, x \geq 0, \, y \geq 0\}$. 
	\end{itemize}
\end{lemma}
The proof of Lemma \ref{Lem-Choice-Alpha} is given in Section \ref{Subsec-Lem-Alpha}.
\\
Let us now fix $\varepsilon >0$ and take $\alpha$ as in Lemma \ref{Lem-Choice-Alpha}. We then define the following oriented paths:
\begin{equation}
	\label{Contours}
	\left\{ 
		\begin{array}{lll}
			\mathcal{C}_1 &=& \{ X_{\alpha,+}(\tilde x), \, \tilde x \hbox{ from } L_0 \hbox{ to } 0 \},
			\\
			\mathcal{C}_2 &=& \{ \overline{X_{\alpha,-}(\tilde x)}, \, \tilde x  \hbox{ from } 0 \hbox{ to } -L_0 \},
			\\
			\mathcal{C}_3 &=& \{ X_{\alpha,-}(\tilde x), \, \tilde x  \hbox{ from } -L_0 \hbox{ to } 0 \},	
			\\
			\mathcal{C}_4 &=& \{ \overline{X_{\alpha,+}(\tilde x)}, \, \tilde x  \hbox{ from } 0 \hbox{ to } L_0 \}, 
		\end{array}
	\right.
	\qquad
	\mathcal{C}  =  \mathcal{C}_1 \cup \mathcal{C}_2 \cup \mathcal{C}_3 \cup \mathcal{C}_4.
\end{equation}
The contour $\mathcal{C}$ is a closed path contained in $\mathcal{S}(L_0(1+\varepsilon)) \setminus \mathcal{S}(L_0)$, see Figure \ref{Fig3}.
Besides, it is easy to check that $\overline{X_{\alpha,\pm}(\tilde x)} = X_{-\alpha,\mp}(\tilde x)$ for all $\tilde x \in [-L_0, L_0]$. This suggests that, to reach functions which are analytic in $\mathcal{S}(L_0(1+\varepsilon))$, one should use the operators $K_\alpha$ and $K_{-\alpha}$. Indeed, using both these operators, we get the following result, proved in Section \ref{Subsec-Proof-Anal}:
\begin{proposition}
	\label{Prop-Analytic-In-Omega}
	Let $L>0$, $L_0>L$, $\varepsilon >0$, and $\alpha$ as in Lemma \ref{Lem-Choice-Alpha}.
	\\
	Then, for any function $k$ defined on $(-L,L)$ which can be extended analytically on $\overline{\mathcal{S}(L_0(1+\varepsilon))}$, one can find two functions $h_+ \in L^2(-L_0,L_0)$ and $h_- \in L^2(-L_0, L_0)$ such that the function $k - K_{\alpha}(h_+) - K_{-\alpha} (h_-)$ can be extended analytically on the ball of radius $L_0$.
\end{proposition}
Combining this result with Proposition \ref{Prop-K-alpha-0}, we get the following immediate corollary:
\begin{corollary}
	\label{Cor-Almost-The-Result}
	Let $L>0$.
	For any $L_0>L$ and $\varepsilon >0$, there exists a continuous function $\alpha: (-L_0, L_0) \to \mathbb{R}$ such that any function $k$ defined on $(-L,L)$ which can be extended analytically on $\overline{\mathcal{S}(L_0(1+\varepsilon))}$ can be decomposed as 
	\begin{equation}
		k(x) = K_{\alpha} (h_+)(x) + K_{-\alpha}(h_-) + K_{0,L_1} (h_0)
		\quad \hbox{�for } x \in (-L,L),
	\end{equation}
	with $L_1 = (L+L_0)/2$, for some $h_+ \in L^2(-L_0,L_0)$, $h_- \in L^2(-L_0, L_0)$, and $h_0 \in L^2(-L_1,L_1)$.
\end{corollary}
Theorem \ref{Thm-Main} is then an immediate consequence of Corollary \ref{Cor-Almost-The-Result} and Theorem \ref{Thm-Reachable-In-Terms-Of-Source-Term}.
\subsection{Proof of Proposition \ref{Prop-K-alpha-0}}\label{Subsec-Proof-K-alpha-0}
\begin{proof}
	We start by writing the operator $K_{0,{L_0}}$ in \eqref{Def-K-0} slightly differently:
	\begin{equation}
		K_{0,{L_0}}(h) (x) = \frac{2}{L_0^2 \sqrt{\pi}} \int_{-{L_0}}^{L_0} \frac{h(\tilde x)}{x^2/L_0^2 - 2 x \tilde x/L_0^2 + 1} d \tilde x.
	\end{equation}
	Therefore, using \eqref{Generatrice-Series}, we get
	\begin{equation}
		\label{K-0-Cheby}
		K_{0,{L_0}}(h) (x) = \frac{2}{L_0^2 \sqrt{\pi}} \sum_{n \geq 0} \left(\frac{x}{{L_0}}\right)^n \int_{-{L_0}}^{L_0} U_n\left( \frac{\tilde x}{{L_0}} \right) h (\tilde x) d \tilde x.
	\end{equation}
	We then recall that the Chebychev polynomials are orthogonal for the scalar product $L^2(\sqrt{1-x^2} \, dx)$, i.e. for all $m$ and $n$ in $\mathbb{N}$, 
	\begin{equation}
		\int_{-1}^1 U_n(\tilde x) U_m(\tilde x) \sqrt{1 - \tilde x^2} \, d\tilde x = \frac{\pi}{2} \delta_{n,m},
	\end{equation}
	where $\delta_{n,m}$ is the Kronecker symbol, so that we have in particular
	\begin{equation}
		\label{Dual-Cheb}
		\int_{-{L_0}}^{L_0} U_n\left( \frac{\tilde x}{{L_0}} \right)U_m\left( \frac{\tilde x}{{L_0}} \right) \sqrt{1 - \frac{\tilde x^2}{L_0^2}} \, d\tilde x = \frac{\pi {L_0} }{2} \delta_{n,m}.
	\end{equation}
	Let us now consider a function $k$ which can be extended analytically on the closure of the ball of radius ${L_0}$. Then $k$ is characterized by its power series expansion:
	$$
		k(z) = \sum_{n \geq 0} k_n z^n, \quad z \in \overline{B(0, {L_0})},
	$$
	and the coefficients $k_n$ satisfy
	\begin{equation}
		\label{Cond-coeff-k-n}
		\sum_{n \geq 0} |k_n| L_0^n < \infty, 
	\end{equation}
	Therefore, using \eqref{K-0-Cheby} and \eqref{Dual-Cheb}, one easily checks that a good candidate $h$ for solving $K_{0,{L_0}}(h) = k$ is given by 
	\begin{equation}
		\label{Def-h-candidate}
		h(\tilde x) 
		= 
		\frac{{L_0}}{\sqrt{\pi}}
		\sum_{m \geq 0} L_0^m k_m U_m\left( \frac{\tilde x}{{L_0}} \right) \sqrt{1 - \frac{\tilde x^2}{L_0^2}}, \quad \tilde x \in (-{L_0}, {L_0}).
	\end{equation}
	We then check that $h$ indeed belongs to $L^2(-{L_0}, {L_0})$. This follows from the following computations, based on \eqref{Dual-Cheb}:
	\begin{eqnarray*}
		\int_{-{L_0}}^{L_0} |h(\tilde x)|^2 \, d\tilde x
		& \leq &
		\int_{-{L_0}}^{L_0} |h(\tilde x)|^2 \, \frac{1}{\sqrt{1 - \tilde x^2/L_0^2}}d\tilde x
		\\
		& = & 
		\frac{L_0^2}{\pi} \int_{-{L_0}}^{L_0} \sum_{m,n \geq 0} L_0^{m+n}  k_m \overline{k_n} U_m\left( \frac{\tilde x}{{L_0}} \right)U_n\left( \frac{\tilde x}{{L_0}} \right) \sqrt{1 - \tilde x^2/L_0^2} \, d\tilde x
		\\�
		& = &
		\frac{{L_0}}{2} \sum_{m \geq 0}  L_0^{2m} |k_m|^2, 
	\end{eqnarray*}
	which is finite due to \eqref{Cond-coeff-k-n}. This concludes the proof of Proposition \ref{Prop-K-alpha-0}.
\end{proof}
\begin{remark}
	Note that the quantity 
	$$
		\sum_{m \geq 0}  L_0^{2m} |k_m|^2
	$$
	appearing in the proof is related to the norm of the function $k$ in the Hardy space $\mathcal{H}^2$ on the ball of radius $L_0$.
\end{remark}
\subsection{Proof of Lemma \ref{Lem-Choice-Alpha}}\label{Subsec-Lem-Alpha}
\begin{proof}
Let us first remark that rescaling if needed as follows
$$
	\alpha(\tilde x)  \longleftrightarrow L_0^2 \tilde \alpha\left( \tau\right), \hbox{ with } \tau = \frac{\tilde x}{L_0},
$$
we can focus on the case $L_0 = 1$ without loss of generality. 
Indeed, in that case, 
$$
	\tilde X_{\tilde \alpha,+}(\tau) =\tau+ \ic \sqrt{1 - \tau^2 + \ic \tilde \alpha \left( \tau \right)}
	= L_0^2 X_{\alpha,+}(\tilde x), \hbox{ with } \tilde x = L_0 \tau.
$$
In the following, we call $\tau$ the rescaled variable and we simply denote the rescaled functions $\tilde \alpha,\, \tilde X_{\tilde \alpha,+}$ by $\alpha, X_{\alpha,+}$ to simplify notations.
\\
For $p \in \mathbb{N} \setminus \left\lbrace 0 \right\rbrace$, we set
\begin{equation}
	\label{Def-alpha-p}
	\alpha_p := \tau\in (-1,1) \mapsto 2\, \vert \tau \vert \left( 1 - \tau^{2 p} \right).
\end{equation}
Note that with this choice, item \emph{(i)} of Lemma \ref{Lem-Choice-Alpha} is obvious.
\\
With this choice, we also immediately get that 
$$
	X_{\alpha_p,-}(\tau) = - X_{\alpha_p,+}(- \tau), \quad \tau \in [-1,1],
$$
i.e. item \emph{(ii)} of Lemma \ref{Lem-Choice-Alpha}.
\\
In order to study the map $X_{\alpha_p, +}$ on $[-1,1]$, it is therefore sufficient to characterize the sets of $\{X_{\alpha_p, +}(\tau), \, \tau \in [0,1]\}$ and $\{X_{\alpha_p, -}(\tau), \, \tau \in [0,1]\} = - \{X_{\alpha_p, -}(\tau), \, \tau \in [-1,0]\} $.
\\
If we define, for $\tau \in (-1,1)$,
$$
	\gamma(\tau)  = \sqrt{(1-\tau^2)^2 + \alpha_p(\tau)^2} 
	\text{ and } 
	\theta(\tau) \in \left[0,\frac{\pi}{2}\right]
	\text{ s.t. } 
	\cos(\theta(\tau)) = \frac{1-\tau^2}{\gamma(\tau)},
	\ 
	\sin(\theta(\tau)) = \frac{\alpha_p(\tau)}{\gamma(\tau)},
$$
which can be extended continuously for $\tau = \pm 1$, we obtain
$$
	X_{\alpha_p,+}(\tau) 
	= 
	\tau - \sqrt{\gamma(\tau)} \sin\left( \frac{\theta(\tau)}{2} \right) + \ic \sqrt{\gamma(\tau)} \cos \left( \frac{\theta(\tau)}{2} \right)
$$
implying in particular that $\Im(X_{\alpha_p,+}(\tau)) \geq 0$ for all $\tau \in (-1,1)$.
\\
Furthermore, using that
$$
	\sin \left( \frac{\theta(\tau)}{2} \right) = \sqrt{\frac{1}{2}( 1 - \cos(\theta(\tau))} = \frac{1}{\sqrt{2 \gamma(\tau)}} \sqrt{\gamma(\tau) - (1- \tau^2))},
$$
we get
\begin{align*}
	\vert X_{\alpha_p,+}(\tau) \vert^2  
	& 
	= 
	\tau^2 - 2\, \tau\, \sqrt{\gamma(\tau)} \sin \left( \frac{\theta(\tau)}{2} \right) + \gamma(\tau) 
	\\
	&
	= 
	\tau^2 + \gamma(\tau) - \sqrt{2} \tau \sqrt{\gamma(\tau) - (1-\tau^2)} 
	\\
	 &  = 1 + \sqrt{\gamma(\tau) - (1-\tau^2)} \left( \sqrt{\gamma(\tau) - (1-\tau^2)} - \sqrt{2} \tau \right).
\end{align*}
Under this form, we clearly have that for $\tau\in (-1,0)$, $\vert X_{\alpha_p,+}(\tau) \vert  > 1$. Besides,  for $\tau\in (0,1)$, we have
$$
	\sqrt{\gamma(\tau) - (1-\tau^2)} <  \sqrt{2} \tau
	 \Leftrightarrow 
	 \gamma(\tau) < 1 + \tau^2 
	\Leftrightarrow 
	\gamma(\tau)^2 < (1+\tau^2)^2 
	\Leftrightarrow 
	\alpha_p(\tau)^2 < 4\, \tau^2,
$$
the last inequality being obviously true for any $\tau \in (0,1)$, recall the definition of $\alpha_p$ in \eqref{Def-alpha-p}. This implies that $\vert X_{\alpha_p,+}(\tau) \vert < 1$ for all $\tau \in (0,1)$.
\\
We can then remark that 
$$
	\Re(X_{\alpha_p,+}(\tau)) =  \tau - \sqrt{\gamma(\tau)} \sin\left( \frac{\theta(\tau)}{2} \right)
	= \tau - \frac{1}{\sqrt{2}} \sqrt{\gamma(\tau) - (1-\tau^2)},
$$
so that similarly as above, $\Re(X_{\alpha_p,+}(\tau))> 0$ for $\tau \in (0,1)$ and $\Re(X_{\alpha_p,+}(\tau))< 0$ for $\tau \in (-1,0)$.
\\
We have thus proved that $\{ X_{\alpha_p,+}(\tau), \, \tau \in (-1,0) \}$ is contained in $\{x +iy, \, x \leq 0, \, y \geq 0\} \setminus B(0,1)$. Using  item \emph{(ii)} of Lemma \ref{Lem-Choice-Alpha}, one easily checks that the set $\{X_{\alpha,-}(\tau), \, \tau \in [0,1]\}$ is included in the set $\{x + \ic y, \, x \geq 0, \, y \leq 0\}\setminus B(0, 1)$.
\\
Therefore, to finish the proof of item \emph{(iii)} of Lemma \ref{Lem-Choice-Alpha}, we only need to prove that $\{ X_{\alpha_p,+}(\tau), \, \tau \in (-1,0) \}$ describes a rectifiable curve. This can be done easily by a tedious computation after having noticed that 
$$
	X_{\alpha_p,+}(\tau) =  \tau - \frac{1}{\sqrt{2}} \sqrt{\gamma(\tau) - (1-\tau^2)} + \ic \frac{1}{\sqrt{2}} \sqrt{\gamma(\tau) + (1-\tau^2)}.
$$
The details of the computations are left to the readers.
\\ \par
We shall now focus on the proof of item \emph{(iv)} of Lemma \ref{Lem-Choice-Alpha}.
\\
As $\Re(X_{\alpha_p,+}(\tau)) \geq 0$ for $\tau \in [0,1]$, we have, for all $\tau \in [0,1]$,
\begin{align*}
	\vert \Re(X_{\alpha_p,+}(\tau)) \vert  + \vert \Im(X_{\alpha_p,+}(\tau)) \vert 
	& 
	= \tau + \sqrt{\gamma(\tau)} \left( \cos\left(\frac{\theta(\tau)}{2}\right) - \sin \left( \frac{\theta(\tau)}{2} \right) \right) 
	\\
	& 
	= 1 + (\tau-1 ) + \frac{1}{\sqrt{2}} \left[ \sqrt{\gamma(\tau) + (1-\tau^2)}  - \sqrt{\gamma(\tau) - (1-\tau^2)} \right].
\end{align*}
For $\tau \in (0,1)$, we hence have
\begin{align*}
	\vert \Re(X_{\alpha_p,+}(\tau)) \vert  + \vert \Im(X_{\alpha_p,+}(\tau)) \vert 
	\geq  1 & \Leftrightarrow \sqrt{\gamma(\tau) + (1-\tau^2)}  - \sqrt{\gamma(\tau) - (1-\tau^2)} \geq \sqrt{2}(1-\tau) 
	\\
	& \Leftrightarrow 2 \gamma(\tau) - 2 \sqrt{\gamma(\tau)^2 - (1-\tau^2)^2} \geq 2 (1-\tau)^2 
	\\
	& \Leftrightarrow \gamma(\tau) - \alpha_p(\tau) \geq  (1-\tau)^2 
	\\
	&  \Leftrightarrow \gamma(\tau)^2 \geq \left[ (1-\tau)^2 + \alpha_p(\tau) \right]^2 
	\\
	& \Leftrightarrow (1-\tau^2)^2 \geq (1-\tau)^4 + 2\, \tau \, (1-\tau^{2\,p})(1-\tau)^2 
	\\
	& \Leftrightarrow 2\,{\left( \tau - 1\right) }^{2}\,\tau\,\left( {\tau}^{2\,p}+1\right)  \geq 0.
\end{align*}
This last inequality obviously holds true, so we have 
\begin{equation}
	\forall \tau \in [0,1], 
	\quad
	\vert \Re(X_{\alpha_p,+}(\tau)) \vert  + \vert \Im(X_{\alpha_p,+}(\tau)) \vert 
	\geq  1.
\end{equation}
Let us then define
$$
	g_p := \tau \in (0,1) 
		\mapsto 
	\vert \Re(X_{\alpha_p,+}(\tau)) \vert  + \vert \Im(X_{\alpha_p,+}(\tau)) \vert - 1 
	= 
	(\tau-1 ) + \frac{1}{\sqrt{2}} \left[ \sqrt{\gamma(\tau) + (1-\tau^2)}  - \sqrt{\gamma(\tau) - (1-\tau^2)} \right].
$$
We already know that for all $\tau \in [0,1]$, $g_p(\tau) \geq 0$. Our next goal is to show that in fact, $g_p$ is bounded on $[0,1]$ by some bounds going to $0$ as $p \to \infty$. In order to do that, we will decompose the interval $[0,1]$ in two intervals $[0, \tau_p]$ and $[\tau_p,1]$ for some parameters $\tau_p \in [0,1]$ going to $1$ as $p \to \infty$, and we will establish bounds going to $0$ as $p \to \infty$ on each of these intervals.
\\
Let us start by the following remark: for $\tau_p > 1/\sqrt{2}$, for all $\tau \in [\tau_p,1)$, 
$$
	g_p(\tau) = (\tau-1) + \frac{1}{\sqrt{2}} \left[ \sqrt{\gamma(\tau) + (1-\tau^2)} - 
\sqrt{\gamma(\tau) - (1-\tau^2)} \right] \leq (\tau-1) + \sqrt{1-\tau^2} \leq \tau_p - 1 + \sqrt{1-\tau_p^2}.
$$
Therefore, if 
\begin{equation}
	\label{Condition-First-Interval}
	\lim_{p \rightarrow \infty} \tau_p = 1,
\end{equation}
we get  $\lim_{p\to \infty} \|g_p \|_{L^\infty(\tau_p,1)} = 0$.
\\
Secondly, for any $\tau \in (0,\tau_p)$, we observe that
$
	\gamma(\tau)^2 = (1+\tau^2)^2 + r(\tau)
$
with 
$
	r(\tau) = -8\, \tau^{2p+2} \left(1- \frac{1}{2} \tau^{2p} ) \right)
$
which implies in particular that
$$
	\sup_{\tau \in [0, \tau_p]} \left\vert \frac{r(\tau)}{\tau^2} \right\vert \leq 8 \tau_p^{2p}.
$$
Therefore, if we choose $\tau_p$ such that 
\begin{equation}
	\label{Cond-SecondInterval}
	\lim_{p \to \infty} \tau_p^{2p} = 0,
\end{equation}
we have, for all $\tau \in [0, \tau_p]$
$$
	\left| 
		\gamma(\tau) - \left( (1+\tau^2) + \frac{r(\tau)}{2(1+\tau^2)}\right)
	\right|
	\leq
	C r(\tau), 
$$
for some $C$ independent of $p$.
\\ 
This leads, still under condition \eqref{Cond-SecondInterval} that for $\tau \in [0,\tau_p]$,
$$
	\left|
		\frac{1}{\sqrt{2}} \left[ \sqrt{\gamma(\tau) + (1-\tau^2)}   - \sqrt{\gamma(\tau) - (1-\tau^2)}\right]
	- \left( 1 - \tau + \frac{1}{8\,(1+\tau^2)} \left[ r(\tau) - 2\frac{r(\tau)}{\tau^2} \right] \right)
	\right| 
	\leq C \frac{r(\tau)}{\tau^2}, 
$$
for some $C$ independent of $p$.
\\
Therefore, under condition \eqref{Cond-SecondInterval}, we have
$$
	\sup_{\tau\in [0,\tau_p]} g_p(\tau) \leq C \sup_{\tau\in[0,\tau_p]} \left\vert \frac{r(\tau)}{\tau^2} \right\vert \underset{p\to \infty}{\rightarrow} 0, 
$$
\par
\begin{figure}[ht]
\begin{center}
\includegraphics[width = 0.7\textwidth,height = 0.3\textheight]{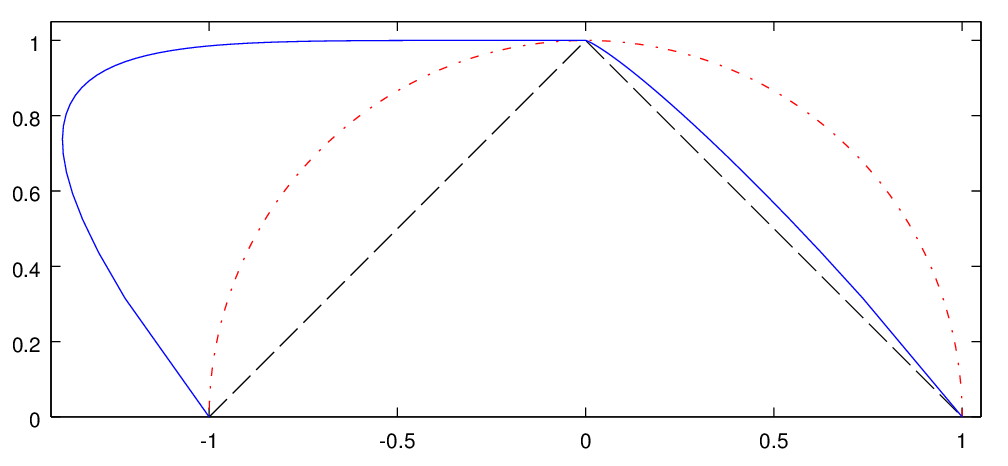}
\end{center}
\caption{In red, $.-$ the Euclidean sphere; in black, $--$ the boundary of the $\ell^1(\mathbb{R}^2)$ ball; in blue, the curve $X_{\alpha_2,+}$ for $\tau \in [-1,1]$.}
\label{Fig1}
\end{figure}
\begin{figure}[ht]
\begin{center}
\includegraphics[width = 0.5\textwidth, height = 0.3\textheight]{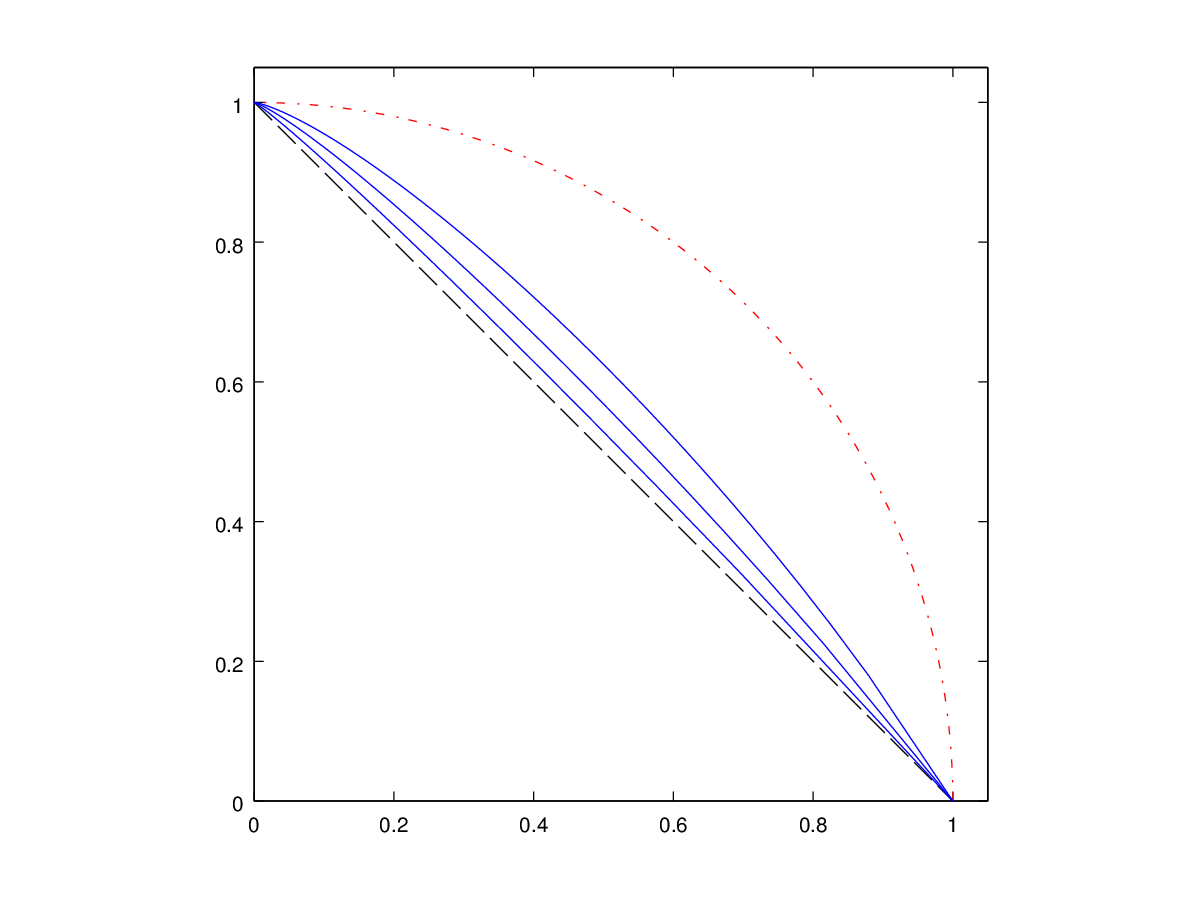}
\end{center}
\caption{$X_{\alpha_p,+}$ for $\tau \in (0,1)$, and $p = 1$, $2$ and $5$.}
\label{Fig2}
\end{figure}

We thus choose $\tau_p = 1 - 1/\sqrt{p}$, so that conditions \eqref{Condition-First-Interval} and \eqref{Cond-SecondInterval} are satisfied, and we obtain
that
$$
	\lim_{p \to \infty} \| g_p \|_{L^\infty(0,1)} = 0.
$$
Therefore, for all $\varepsilon >0$, we can choose $p \in \mathbb{N}$ such that $ \| g_p \|_{L^\infty(0,1)} \leq \varepsilon$. This means geometrically that the set $\{X_{\alpha,+}(\tau), \, \tau \in [0,1]\}$ is included in $\mathcal{S}(1+\varepsilon) \setminus \mathcal{S}(1)$.
\\
Finally, the fact that $\{X_{\alpha,+}(\tau), \, \tau \in [0,1]\}$ is a rectifiable curve can be done as in the proof of item \emph{(iii)} of Lemma \ref{Lem-Choice-Alpha} by explicit computations. This finishes the proof of item \emph{(iv)} of Lemma \ref{Lem-Choice-Alpha}.
\\
The proof of Lemma \ref{Lem-Choice-Alpha} is now completed. We finish it with Figures \ref{Fig1}--\ref{Fig2} illustrating Lemma \ref{Lem-Choice-Alpha}.
\end{proof}
\subsection{Proof of Proposition \ref{Prop-Analytic-In-Omega}}\label{Subsec-Proof-Anal}
Let $L>0$, $L_0 > L$, $\varepsilon >0$, and $\alpha$ as in Lemma \ref{Lem-Choice-Alpha}. 
\\
Let us then consider a function $k$ which can be extended analytically on $\overline{\mathcal{S}(L_0(1+\varepsilon))}$. We still denote by $k$ its analytic expansion.
\\
We note that the oriented path $\mathcal{C}$ in \eqref{Contours} is included in ${\mathcal{S}(L_0(1+\varepsilon))}$. We can therefore use Cauchy's integral formula:	
\begin{equation}	
	\forall x \in [-L,L], \quad 
	k(x) = \frac{1}{2\, \ic \pi }\int_\mathcal{C} \frac{k(z)}{z-x}\, dz, 
\end{equation}
which in our context yields:
\begin{eqnarray}
	k(x) & = \ds \frac{1}{2\, \ic \pi} 
	&
	\left[ 
		- 
		\int_0^{L_0} \frac{k(X_{\alpha,+}(\tilde x))}{X_{\alpha,+}(\tilde x)  - x} X_{\alpha,+}'(\tilde x) \, d\tilde x 
		+ 
		\int_{-L_0}^0\frac{k(X_{\alpha-}(\tilde x))}{X_{\alpha,-}(\tilde x)  - x} X_{\alpha,-}'(\tilde x) \, d\tilde x 
	\right. 
	\label{Cauchy-k}
	\\
	& & \left.  
		+ 
		\int_0^{L_0} \frac{k\left(\overline{X_{\alpha,+}(\tilde x)}\right)}{\overline{X_{\alpha,+}(\tilde x)} - x} \overline{X_{\alpha,+}'(\tilde x)} d\, \tilde x 
		- 
		\int_{-L_0}^0 \frac{k\left(\overline{X_{\alpha,-}(\tilde x)}\right)}{\overline{X_{\alpha,-}(\tilde x)} - x}  \overline{X_{\alpha,-}'(\tilde x)} d\,\tilde x  	\right].
	\notag
\end{eqnarray}
\par
Let us then recall that for $h_+ \in L^2(-L_0, L_0)$ and $h_- \in L^2(-L_0, L_0)$, we have
\begin{eqnarray}	
	K_{\alpha} (h_+) (x)
		& = &
		\frac{2}{\sqrt{\pi}}
		\int_{-L_0}^{L_0}
		\frac{h_+(\tilde x)}{X_{\alpha,+}(\tilde x) -X_{\alpha,-}(\tilde x) } 
		\left(
			\frac{1}{x- X_{\alpha,+}(\tilde x)} 
			- 
			 \frac{1}{x- X_{\alpha,-}(\tilde x)}
		\right)
		d \tilde x
		\label{K-alpha-+}
	\\
	K_{-\alpha} (h_-) (x)
		& = &
		\frac{2}{\sqrt{\pi}}
		\int_{-L_0}^{L_0}
		\frac{h_-(\tilde x)}{X_{-\alpha,+}(\tilde x) -X_{-\alpha,-}(\tilde x) } 
		\left(
			\frac{1}{x- X_{-\alpha,+}(\tilde x)} 
			- 
			 \frac{1}{x- X_{-\alpha,-}(\tilde x)}
		\right)
		d \tilde x	
		\notag
	\\
		& = &
		\frac{2}{\sqrt{\pi}}
		\int_{-L_0}^{L_0}
		\frac{h_-(\tilde x)}{\overline{X_{\alpha,-}(\tilde x)} -\overline{X_{\alpha,+}(\tilde x) }} 
		\left(
			\frac{1}{x- \overline{X_{\alpha,-}(\tilde x)}} 
			- 
			 \frac{1}{x- \overline{X_{\alpha,+}(\tilde x) }}
		\right)
		d \tilde x.	
		\label{K-alpha--}
\end{eqnarray}
In view of \eqref{Cauchy-k}, it is therefore natural to choose $h_+$ such that 
\begin{equation}
	\label{Def-h-+}
	\frac{2}{\sqrt{\pi}}\frac{h_+(\tilde x)}{X_{\alpha,+}(\tilde x) -X_{\alpha,-}(\tilde x) } 
	= 
	\left\{
		\begin{array}{ll}
			\ds 
			\frac{1}{2 \ic \pi} k(X_{\alpha,+}(\tilde x))X_{\alpha,+}'(\tilde x) &\hbox{ for } \tilde x > 0, 
			\smallskip
			\\
			\ds - \frac{1}{2 \ic \pi} k(X_{\alpha,-}(\tilde x))X_{\alpha,-}'(\tilde x) &\hbox{ for } \tilde x \leq 0, 
		\end{array}
	\right.
\end{equation}
and $h_-$ such that
\begin{equation}
	\label{Def-h-}
	\frac{2}{\sqrt{\pi}} \frac{h_-(\tilde x)}{\overline{X_{\alpha,-}(\tilde x)} -\overline{X_{\alpha,+}(\tilde x) }}
	= 
		\left\{
		\begin{array}{ll}
			\ds 
			- \frac{1}{2 \ic \pi} k\left(\overline{X_{\alpha,+}(\tilde x)}\right)\overline{X_{\alpha,+}'(\tilde x)}
			& \hbox{ for } \tilde x > 0, 
			\smallskip
			\\
			\ds  
			\frac{1}{2 \ic \pi}k\left(\overline{X_{\alpha,-}(\tilde x)}\right) \overline{X_{\alpha,-}'(\tilde x)}
			& \hbox{ for } \tilde x \leq 0. 
		\end{array}
	\right.
\end{equation}
\par
Let us then check that the two above definitions \eqref{Def-h-+}--\eqref{Def-h-} give functions $h_+$, $h_-$ in $L^2(-L_0, L_0)$. We explain in details how to show that $h_+ \in L^2(0,L_0)$.
\\
On $(0,L_0)$, we have
$$
	h_+ (\tilde x ) = \frac{1}{4 \ic \sqrt{\pi} } \left(X_{\alpha,+}(\tilde x) -X_{\alpha,-}(\tilde x)\right) k(X_{\alpha,+} (\tilde x)) X_{\alpha,+}'(\tilde x).
$$
But for all $\tilde x \in (0,L_0)$, $X_{\alpha,+}(\tilde x) \in \mathcal{S}(0,L_0(1+\varepsilon))$ on which $k$ is bounded (since it is analytic on $\overline{\mathcal{S}(L_0(1+\varepsilon))}$). Therefore, we only have to check that 
$$
	(X_{\alpha,+}(\tilde x) -X_{\alpha,-}(\tilde x )) X_{\alpha,+}'(\tilde x)
$$
belongs to $L^2(0,L_0)$. Explicit computations yield that 
\begin{align*}
	& X_{\alpha,+}'(\tilde x) = 1 + \frac{\ic}{2\sqrt{L^2- \tilde x^2 + \ic \alpha(\tilde x)}} \left(- 2 \tilde x + \alpha'(\tilde x) \right) ,
	\\
	& X_{\alpha,+}(\tilde x) -X_{\alpha,-}(\tilde x ) = 2 \ic \sqrt{L^2- \tilde x^2 + \ic \alpha(\tilde x)}, 
\end{align*}
so that for all $\tilde x \in (0,L_0)$,
$$
	(X_{\alpha,+}(\tilde x) -X_{\alpha,-}(\tilde x )) X_{\alpha,+}'(\tilde x)
	= 
	2 \ic \sqrt{L^2- \tilde x^2 + \ic \alpha(\tilde x)} -  \left(- 2 \tilde x + \alpha'(\tilde x) \right), 
$$
which is obviously bounded in view of item \emph{(i)} of Lemma \ref{Lem-Choice-Alpha}. Therefore, $h_+$ given by \eqref{Def-h-+} belongs to $L^2(0, L_0)$.
\\
Of course, similar computations can be done to show that $h_+ \in L^2(-L_0, 0)$ and $h_- \in L^2(-L_0, L_0)$. The details of these proofs are left to the reader.
\\ \par
Let us then show that the function $k_r$ defined for $x \in (-L,L)$ by 
\begin{equation}
	\label{Def-k-r}
	k_r (x) = k(x) - K_\alpha (h_+)(x) - K_\alpha(h_-)(x)
\end{equation}
can be extended analytically on the ball of size $L_0$. Indeed, from formulae \eqref{Cauchy-k}--\eqref{K-alpha-+}--\eqref{K-alpha--}--\eqref{Def-h-+}--\eqref{Def-h-+}, we have for all $x \in (-L, L)$, 
\begin{align*}
	k_r (x ) = 
		&
		- \frac{2}{\sqrt{\pi}}
		\int_{-L_0}^{L_0}
		\frac{h_+(\tilde x)}{X_{\alpha,+}(\tilde x) -X_{\alpha,-}(\tilde x) } 
		\left(
			\frac{{\bf 1}_{\tilde x <0}}{x- X_{\alpha,+}(\tilde x)} 
			- 
			 \frac{{\bf 1}_{\tilde x >0}}{x- X_{\alpha,-}(\tilde x)}
		\right)
		d \tilde x
		\\
		&
		- \frac{2}{\sqrt{\pi}}
		\int_{-L_0}^{L_0}
		\frac{h_-(\tilde x)}{\overline{X_{\alpha,-}(\tilde x)} -\overline{X_{\alpha,+}(\tilde x) }} 
		\left(
			\frac{{\bf 1}_{\tilde x>0}}{x- \overline{X_{\alpha,-}(\tilde x)}} 
			- 
			 \frac{{\bf 1}_{\tilde x<0}}{x- \overline{X_{\alpha,+}(\tilde x) }}
		\right)
		d \tilde x.	
\end{align*}
But according to Lemma \ref{Lem-Choice-Alpha}, for $\tilde x < 0$, $X_{\alpha,+}(\tilde x) \notin B(0, L_0)$, and for $\tilde x >0$, $X_{\alpha,-}(\tilde x) \notin B(0,L_0)$. Therefore, the singularities in each kernel lie outside $B(0,L_0)$. Therefore, $k_r$ can be extended analytically in $B(0,L_0)$ with the following formula, valid for any $z \in B(0,L_0)$,
\begin{align*}
	k_r (z ) = 
		&
		- \frac{2}{\sqrt{\pi}}
		\int_{-L_0}^{L_0}
		\frac{h_+(\tilde x)}{X_{\alpha,+}(\tilde x) -X_{\alpha,-}(\tilde x) } 
		\left(
			\frac{{\bf 1}_{\tilde x <0}}{z- X_{\alpha,+}(\tilde x)} 
			- 
			 \frac{{\bf 1}_{\tilde x >0}}{z- X_{\alpha,-}(\tilde x)}
		\right)
		d \tilde x
		\\
		&
		- \frac{2}{\sqrt{\pi}}
		\int_{-L_0}^{L_0}
		\frac{h_-(\tilde x)}{\overline{X_{\alpha,-}(\tilde x)} -\overline{X_{\alpha,+}(\tilde x) }} 
		\left(
			\frac{{\bf 1}_{\tilde x>0}}{z- \overline{X_{\alpha,-}(\tilde x)}} 
			- 
			 \frac{{\bf 1}_{\tilde x<0}}{z- \overline{X_{\alpha,+}(\tilde x) }}
		\right)
		d \tilde x.	
\end{align*}
This completes the proof of Proposition \ref{Prop-Analytic-In-Omega}.
\section{Further comments}\label{Sec-Further}
\subsection{The reachable set when the control acts from one side}
One may ask if it is possible to characterize the reachable set of the one-dimensional heat equation controlled from one side only.
\\
To fix the ideas, let $L,\, T >0$ and consider the equation
\begin{equation}
	\label{Eq-Heat-1d-0-L}
		\left\{
			\begin{array}{ll}
				\partial_t u - \partial_{xx} u = 0 & \hbox{ in } (0,T) \times (0,L), 
				\\
				u(t,0) = 0 & \hbox{ in } (0,T), 
				\\
				u(t,L) = v(t)  & \hbox{ in } (0,T), 
				\\
				u(0, x) = 0 & \hbox{ in } (0,L).
			\end{array}
		\right.
\end{equation}
In this context, we define the reachable set $\mathscr{R}_{L, u(t,0) = 0}(T)$ at time $T>0$  as follows:
\begin{equation}
	\label{Def-Reachable-T-Onesidecontrol}
		\mathscr{R}_{L, u(t,0) = 0} (T)
		=
		\{u(T)\ |\  u \hbox{ solving \eqref{Eq-Heat-1d-0-L} with control functions } v \in L^2(0,T) \}.
\end{equation}
Again, this set is a vector space independent of the time $T>0$ and we therefore simply write $\mathscr{R}_{L, u(t,0) = 0}$ instead of $\mathscr{R}_{L, u(t,0) = 0} (T)$.
\\
As a corollary of Theorem \ref{Thm-Main}, one can prove the following result:
\begin{theorem}
	\label{Thm-OneSideControl}
	Any function $u \in L^2(0,L)$ whose odd extension to $(-L,L)$ can be extended analytically to $\overline{\mathcal{S}(L)}$ belongs to $\mathscr{R}_{L, u(t,0) = 0}$.
\end{theorem}
Theorem \ref{Thm-OneSideControl} is in fact an immediate consequence of Theorem \ref{Thm-Main}. Indeed, if $u_1 \in L^2(0,L)$ has an odd extension $\tilde u_1$ to $(-L,L)$ which can be extended analytically to $\overline{\mathcal{S}(L)}$, then $\tilde u_1 \in \mathscr{R}(L)$ from Theorem \ref{Thm-Main}. If we denote by $\tilde u$ a corresponding trajectory of \eqref{Eq-Heat-1d} starting from $\tilde u(0, \cdot) = 0$ in $(-L,L)$, taking value $\tilde u_1$ at time $T$ in $(-L,L)$ and having control functions $v_-, \, v_+ \in L^2(0,T)$, one can check that for $(t,x) \in (0,T) \times (0,L)$, $u (t,x) = (\tilde u(t,x) - \tilde u(t,-x))/2$ solves \eqref{Eq-Heat-1d-0-L} with control function $v(t) = (v_+ (t)- v_-(t))/2$ and its value at time $T$ is $u_1$ in $(0,L)$, i.e. $u_1 \in \mathscr{R}_{L, u(t,0) = 0}$.
\\
Also note that Theorem \ref{Thm-OneSideControl} is mainly sharp as \cite[Theorem 1]{MartinRosierRouchon-2016} states that any state in $\mathscr{R}_{L, u(t,0) = 0}$ should have an odd extension which can be extended analytically to the set $\mathcal{S}(L)$.
\subsection{The multi-dimensional case}
The Carleman estimate stated in Theorem \ref{Thm-Carleman} can be easily generalized to heat equations in spatial domains $\Omega$ which are multi-dimensional Euclidean balls, and with observation on the whole sphere. However, it is not clear how to use it in a clever way to get a sharp description of the reachable set. This issue will be studied in a forthcoming work, as well as simple geometries like strips.

%
%
%
\bibliographystyle{plain}

\end{document}